\documentclass[11pt,reqno]{amsart}      
\usepackage{amssymb,amsthm}             
\numberwithin{equation}{section}        

\renewcommand{\Re}{{\mathbb R}}         
\newcommand{\la}{\langle}               
\newcommand{\ra}{\rangle}               
\newcommand{\half}{\frac{1}{2}}         
\newcommand{\eps}{\epsilon}
\newcommand{\Ric}{\text{\rm Ric}}	
\newcommand{\Ham}{\mathcal H}		
\newcommand{\tr}{\text{\rm tr}}		
\newcommand{\semidir}{\ltimes}          
\newcommand{\dev}{\mathcal D}		
\newcommand{\FF}{\mathcal F}		
\newcommand{\UU}{\mathcal U}		
\newcommand{\Vol}{\text{\rm Vol}}	
\newcommand{\hK}{\hat K}		
\newcommand{\length}{\text{\em length}}
\newcommand{\Lie}{\text{\rm Lie}}
\newcommand{\Ad}{\text{\rm Ad}}
\newcommand{\vect}{\text{\rm vec}}
\newcommand{\genus}{\text{\rm genus}}

\newcommand{\ISO}{\text{\rm ISO}}
\newcommand{\Orth}{\mathrm O}
\newcommand{\SO}{\text{\rm SO}}
\newcommand{\Sim}{\text{\rm Sim}}
\newcommand{\bLGam}{\overline{L(\Gamma)}}
\newcommand{\iso}{\mathfrak{iso}}
\newcommand{\so}{\mathfrak{so}}
\newcommand{\frakm}{\mathfrak{m}}
\newcommand{\HH}{\mathbb H}		
\newcommand{\EE}{\mathbb E}
\newcommand{\MM}{\mathcal M}
\newcommand{\CC}{\mathcal C}
\newcommand{\tV}{\tilde V}		
\newcommand{\tM}{\tilde M}
\newcommand{\tN}{\tilde N}
\newcommand{\norm}{\mathbf n}	
\newcommand{\tnorm}{\tilde{\mathbf n}}	
\newcommand{\tvarphi}{\tilde \varphi}	
\newcommand{\id}{\text{\bf Id}}

\theoremstyle{plain}
\newtheorem{thm}{Theorem}[section]
\newtheorem{conj}{Conjecture}
\newtheorem{cor}[thm]{Corollary}
\newtheorem{lemma}[thm]{Lemma}

\newtheorem{definition}[thm]{Definition}
\newtheorem{prop}[thm]{Proposition}

\title[Constant mean curvature foliations]{Constant mean curvature foliations 
of flat space--times} 
\author[L. Andersson]{Lars Andersson}
\thanks{Supported in part by the Swedish Natural
Sciences Research Council (SNSRC),  contract no.  R-RA 4873-307 and NSF,
contract no. DMS 0104402.}
\address{Department of Mathematics\\
University of Miami\\
Coral Gables, FL 33124\\
USA}

\email{larsa\char'100math.miami.edu}

\date{\today \ {\em File: \jobname{.tex}}}

\allowdisplaybreaks[2]

\begin{document}
\begin{abstract}
Let $V$ be a maximal globally hyperbolic flat $n+1$--dimensional 
space--time with compact 
Cauchy surface of
hyperbolic type. We prove that $V$ is globally foliated by constant mean
curvature hypersurfaces $M_{\tau}$, with mean curvature $\tau$ taking all
values in $(-\infty, 0)$. For $n \geq 3$, define the 
rescaled volume of $M_{\tau}$ by $\Ham = |\tau|^n \Vol(M,g)$, where $g$ is
the induced metric. Then $\Ham \geq n^n \Vol(M,g_0)$ where $g_0$ is the
hyperbolic metric on $M$ with sectional curvature $-1$. Equality holds if and
only if $(M,g)$ is isometric to $(M,g_0)$.  
\end{abstract}
\maketitle


\section{Introduction} 

\subsection{Scalar curvature and volume}
Let $M$ be a compact manifold of dimension $n$ 
of hyperbolic type, with 
hyperbolic metric $g_0$ of sectional curvature $-1$. Besson, Courtois and
Gallot, see \cite{gallot:decreas,BCG1,BCG2} and references therein, 
proved that if $g$ is another metric on $M$ with $\Ric_g \geq
\Ric_{g_0}$, then $\Vol(M,g) \geq \Vol(M,g_0)$ with equality if and only if
$g$ is isometric to $g_0$.  In \cite{BCG1} 
the question is also raised whether a
similar conclusion holds under assumptions on the much softer scalar
curvature, and it was proved that the scale invariant quantity 
$$
\int_M |R_g|^{n/2} \mu_g , 
$$
considered as a function on the space of Riemann metrics, 
has a local minimum at the hyperbolic structure.

Let $\CC$ be the space of conformal classes of metrics on a compact manifold 
$M$ of dimension $n$.
The Yamabe invariant of the conformal class $[g] \in \CC$ is 
$$
\mu[g] = \inf_{h \in [g]} \Vol(M,h)^{\frac{2-n}{n}} \int_M R_h \mu_h . 
$$
The $\sigma$--constant of $M$ is defined as 
$$
\sigma(M) = \sup_{[g] \in \CC} \mu[g] . 
$$
In the case we are considering, $M$ does not admit a metric of nonnegative 
scalar curvature, and there is a
unique representative $h \in [g]$ with $R_h = -n(n-1)$. Thus, 
letting $\MM_{-n(n-1)}$ denote the space of metrics with scalar curvature
$-n(n-1)$, 
we have the relations (cf. \cite[Lemma 1]{lebrun:kodaira} for the first
equality)
\begin{align*}
|\sigma(M)|^{n/2} &= \inf_g \int_M |R_g|^{n/2}\mu_g \\
&=  [n(n-1)]^{n/2} \inf_{g \in \MM_{-n(n-1)}} \Vol(M,g)
\end{align*}
This gives a dual formulation of the $\sigma$--constant in terms of volume. 

There is an ambitious program to solve the geometrization problem in the
negative Yamabe class, and as a consequence to 
understand $\sigma(M)$ in dimension $3$, due to M. Anderson,
cf. \cite{anderson:degI} and references therein. The conjecture that one
would like to prove concerning the $\sigma$--constant in the case of
manifolds of hyperbolic type is the following.
\begin{conj}
\label{conj:sigma}
Suppose $M$ is a compact 3--dimensional 
manifold  of hyperbolic type. Then 
the $\sigma$--constant is realized by the hyperbolic metric, i.e. 
\begin{equation}\label{eq:sigmaconj}
\sigma(M) = -6 \Vol(M,g_0)^{2/3}
\end{equation}
where $g_0$ is the hyperbolic metric on $M$ with sectional curvature $-1$. 
\end{conj}
See \cite{anderson:scalar} for the general formulation. 
At present, little is known about $\sigma(M)$ in dimension $3$.
The only case where we have some additional information is in dimension 4. 
Let $M$ be a compact 4--dimensional manifold of hyperbolic type, 
and let $\CC_{flat}$ denote the space of
flat conformal structures. It follows from the Gauss--Bonnet formula that 
\cite[Theorem 8.3]{johnson:millson}
$$
\inf_{g \in \MM_{-12} \cap \CC_{flat}} \Vol(M,g) = \Vol(M,g_0)
$$
The corresponding question appears to be open in dimension 3.

Using Seiberg--Witten theory, LeBrun has proved the following. 
Let $M$ be the underlying manifold of a 
complex surface 
with even first Betti number. Then
\cite{lebrun:kodaira}
$$
|\sigma(M)|^{2} = 32\pi^2 c_1^2(X)
$$
where $X$ is the minimal model of $M$. 

\subsection{The reduced Hamiltonian and volume}\label{sec:redhamvol}
Let $(g,K)$ be 
a Riemann
metric on $g$ on $M$ and a symmetric $2$--tensor $K$. We call $(g,K)$ vacuum
data on $M$ if the vacuum Einstein constraint equations
\begin{subequations}\label{eq:constraint-intro}
\begin{align}
R_g - |K|_g^2 + (\tr_g K)^2 &= 0 \label{eq:constraint-ham-intro} \\
\nabla_a \tr_g K - \nabla^b K_{ab}  &= 0 \label{eq:constraint-mom-intro}
\end{align}
\end{subequations}
are satisfied. Suppose that $\tr_g K$ is constant on $M$. 
Then the quantity 
\begin{equation}\label{eq:Hamdef}
\Ham = |\tr_g K|^n \Vol(M,g)
\end{equation}
is scale free, if we give $K$ the dimension of $(\length)$, as is natural
from the point of view of the Einstein equations, cf. the discussion in
\cite{fischer:moncrief:sigma}. $\Ham$, which may be thought of as a function
$\Ham(g,K)$ on 
the space of solutions to the constraint equations (\ref{eq:constraint-intro}),
is the reduced Hamiltonian for the vacuum 
Einstein equations, in an appropriate gauge. The infimum of $\Ham$ can be
expressed in terms of $\sigma(M)$ in case $M$ does not admit a metric of
nonnegative scalar curvature, 
$$
\inf_{(g,K)}\frac{n-1}{n} \Ham(g,K)^{2/n} = 
- \sigma(M) , 
$$
where the infimum is taken over $(g,K)$ solving (\ref{eq:constraint-intro}). 
We see that $\sigma(M)$ can be studied in terms of the 
volume for $h \in \MM_{-n(n-1)}$ as well as the 
Hamiltonian $\Ham$. 

Questions about $\inf \Ham$ are related in spirit to the positive mass
problem in general relativity, which concerns the Hamiltonian for the vacuum
Einstein equations, in case $M$ is asymptotically Euclidean. In that case,
the mass of $M$ is defined by 
$$
16\pi m(g) = \lim_{r \to \infty} \int_{S_r} (g_{ii,j} - g_{ij,i} )dS^j .
$$
The positive mass theorem states that for AE metrics with nonnegative scalar
curvature, $m(g) \geq 0$, with equality if and only if $g$ is flat. In case
$M$ is compact, the mass in the above sense is zero. 

The $\sigma$--constant conjecture \ref{conj:sigma}, 
leads to the following conjecture, which may be though of as an analogoue to
positive mass for spatially compact space--times.
\begin{conj} Suppose $M$ is a compact 3--dimensional manifold of hyperbolic 
type and let $(g,K)$ be vacuum data on $M$, with constant mean curvature
$\tr K = \tau$. Then the rescaled volume $\Ham(g,K)$ satisfies
\begin{equation}\label{eq:Hamconj}
\Ham(g,K) \geq 3^3 \Vol(M,g_0) ,
\end{equation}
where $g_0$ is the hyperbolic metric on $M$ with sectional curvature
$-1$. Equality holds in (\ref{eq:Hamconj}) if and only if $g$ is isometric to 
$\frac{n^2}{\tau^2} g_0$, and $K = \frac{\tr K}{n} g$. 
\end{conj}
The Hamiltonian $\Ham$ has the property that it is monotone decreasing along
the flow of the Einstein evolution equations, in constant mean curvature
gauge. This was noticed first by Rendall, and studied in detail by Fischer
and Moncrief \cite{fischer:moncrief:sigma} for dimension $3$. 

Here we will study the rescaled volume $\Ham$ in a restricted geometric
setting, namely that of flat space--times. Suppose $M$ is compact of
hyperbolic type. We call $(g,K)$ flat data on $M$,
if 
\begin{subequations}\label{eq:constraint-flat-intro}
\begin{align}
R_{ab} - K_{ac}K^c_{\ b} + \tr K K_{ab} &= 0 \label{eq:gauss} \\
\nabla_a K_{bc} - \nabla_b K_{ac}  &= 0 \label{eq:codazzi}
\end{align}
\end{subequations}
These are just the Gauss and Codazzi equations for a hypersurface in a
Lorentzian flat ambient space. We may view the flat data on $M$ with $\tr K = 
\tau$ for a given $\tau < 0$ as a parametrization of the space of maximal
globally hyperbolic flat
Lorentz structures of topology $M \times \Re$. The 
dimension of the space of flat data, modulo diffeomorphisms, with $\tr K =
\tau$, for a fixed $\tau < 0$, is equal to the dimension of the Zariski
tangent space of the flat conformal structures on $M$.
\begin{thm}\label{thm:Hambound-intro}
Let $M$ be a compact manifold of hyperbolic type, of dimension $n \geq 3$. 
Let $(g,K)$ be flat data on $M$ with $\tr K = \tau$ constant on
$M$, and let $\Ham(g,K)$ be defined by (\ref{eq:Hamdef}). 
Then 
\begin{equation}\label{eq:Hamineq-intro}
\Ham(g,K) \geq n^n \Vol(M,g_0)
\end{equation}
where $g_0$ is the hyperbolic metric on $M$ with sectional curvature
$-1$. Equality holds in (\ref{eq:Hamineq-intro}) if and only if $g$ is isometric to 
$\frac{n^2}{\tau^2} g_0$ and $K = \frac{\tr
K}{n} g$.
\end{thm}
The proof of Theorem \ref{thm:Hambound-intro} is based on the
above--mentioned monotonicity of $\Ham$, together with two
results concerning constant mean curvature foliations of flat space--times. 
The first, cf. section \ref{sec:cmcglobal},
proves the existence of global constant mean curvature foliations. 
\begin{thm}\label{thm:cmcglobal-intro}
Let $V$ be a maximal globally hyperbolic flat space--time with compact Cauchy 
surface $M$ of hyperbolic type. Then $V$ is globally foliated by constant
mean curvature hypersurfaces
\end{thm}
For $(g,K)$ induced data on a constant mean curvature hypersurface of mean
curvature $\tau < 0$, we write $\Ham(\tau) = \Ham(g,K)$. 
\begin{thm}\label{thm:Hamlimit-intro}
Let $V$ be a maximal globally hyperbolic flat space--time with compact Cauchy 
surface of hyperbolic type. Then 
$$
\lim_{\tau \nearrow 0} \Ham(\tau) = 
n^n \Vol(M,g_0)
$$
\end{thm}
Theorem \ref{thm:Hamlimit-intro} together with the above mentioned
monotonicity property proves Theorem \ref{thm:Hambound-intro}. 

\bigskip

\noindent{\bf Acknowledgements:} 
I am grateful to Misha Gromov, Vince Moncrief, Mike Anderson and
Kevin Scannell for helpful discussions. The paper finished during visits to
the Courant institute, New York and Institute Elie Cartan, Nancy. The
hospitality and support during these visits is gratefully acknowledged.

\section{Preliminaries}\label{sec:prel} A $(G,X)$ structure on $V$ is a
marked locally homogenous structure, modelled on a homogenous space $X =
G/H$.
A $(G,X)$--structure on $V$
gives by a standard construction a $(G,X)$--structure on the universal cover
$\tV$, a developing map $\dev: \tV \to X$ and a holonomy homomorphism $\rho:
\pi_1 V \to G$ satisfying the equivariance condition $\dev(\gamma \cdot x) =
\rho(\gamma) \cdot \dev(x)$ for all $\gamma \in \pi_1 V $, $x \in \tV$. 

The $n+1$--dimensional Minkowski space $\Re^{n+1}_1$, is $\Re^{n+1}$ with the 
flat Lorentz metric $-dt^2 + (dx^1)^2 + \cdots + (dx^n)^2$. We will use
coordinates $t,x^1, \dots x^n$ or $x^0, x^1 , x^n$ on $\Re^{n+1}_1$, and
write $|x| = (\sum_{i=1}^n (x^i)^2)^{1/2}$. 
The group of
linear isometries of $\Re^{n+1}_1$ is the Lorentz group $\Orth(n,1)$, and the
group of isometries of $\Re^{n+1}_1$ is the semidirect product $\ISO(n,1) =
\Orth(n,1) \semidir \Re^{n+1}$ of the Lorentz group with the group of
translations of $\Re^{n+1}$. We will denote  by $\Orth(n,1)_+$ the
subgroup of time orientation preserving Lorentz transformations, by
$\SO(n,1)$ the group of Lorentz transformations of determinant $1$, and by 
$\SO(n,1)_0$ the connected
component of the identity in $\Orth(n,1)$. Similarly $\ISO(n,1)_0$ denotes
the connected component of the identity in $\ISO(n,1)$. 
The frame bundle of $\Re^{n+1}_1$
has a natural $\SO(n,1)_0$ reduction. 

Examples of $(G,X)$--structures in dimension $n$  
are 
hyperbolic structures,   with $(G = SO(n,1) , X = \HH^n)$, 
flat conformal structures, with  $(G = \SO(n+1,1), X = S^n)$, 
and flat Lorentz structures, with  $(G = \ISO(n,1), X = \Re^{n+1}_1)$.

\subsection{Deformation spaces}
The space of $(G,X)$ structures on $V$ 
is isomorphic to the deformation space of holonomy representations $\rho:
\pi_1 V \to G$. The formal (Zariski) 
tangent space of the deformation space at $\rho$
is, letting $\Gamma = \rho  \pi_1 V$, given by the cohomology
$H^1(\Gamma, \mathfrak g_{\Ad})$, where $\mathfrak g_{Ad}$ denotes the adjoint action
of $G$ on $\mathfrak g = \Lie(G)$. In general, the deformation space is not a 
manifold, for example the space of flat conformal structures may have 
singularities at the hyperbolic strucure on a compact manifold of hyperbolic
type in dimension $n \geq 4$ \cite{johnson:millson}. Obstructions to integrating a direction in the
Zariski tangent space can be computed in terms of Massey products
\cite{goldman:millson}.

Let $\iso(n,1) = \Lie(\ISO(n,1))$, $\so(n+1,1) = \Lie(\SO(n+1,1))$. 
Note that $\so(n+1,1) = \so(n,1) \oplus \frakm$, $\frakm = \Re^{n+1}$ as a
vector space, and $\iso(n,1) = \so(n,1)\oplus \Re^{n+1}$, where $\Re^{n+1}$
has the trivial Lie algebra structure. The adjoint action of $\SO(n,1)
\subset \SO(n+1,1)$ on $\so(n+1,1)$ induces an action on $M$ which is
equivalent to the vector action of $\SO(n,1)$ on $\Re^{n+1}$. We let
$\Re^{n+1}_{\vect}$ denote the vector representation of $\SO(n,1)$ on
$\Re^{n+1}$. 

Let $\Gamma$ be the holonomy representation for a hyperbolic structure on
$M$. Then 
\begin{align*}
H^1(\Gamma,\so(n+1,1)_{\Ad} ) &= H^1(\Gamma, \so(n,1)_{\Ad} ) 
\oplus H^1(\Gamma , \frakm_{\Ad}), \\
H^1(\Gamma,\iso(n,1)_{\Ad} ) &= H^1(\Gamma, \so(n,1)_{\Ad} ) 
\oplus H^1(\Gamma , \Re^{n+1}_{\vect}) .
\end{align*}
For $n\geq 3$, $H^1(\Gamma, \so(n,1)_{\Ad}) = \{0\}$, by Calabi-Weil
rigidity, so $H^1(\Gamma, \Re^{n+1}_{\vect})$ is the Zariski tangent space for 
both flat conformal structures and flat Lorentz structures. For $n=2$ on the
other hand, the adjoint representation of $\SO(2,1)$ on $\so(2,1)$ is
equivalent to the vector representation on $\Re^3$, so 
$$
H^1(\Gamma, \iso(2,1)_{\Ad}) = H^1(\Gamma, \Re^3_{\vect}) \oplus  H^1(\Gamma,
\Re^3_{\vect}) . 
$$
Further, in case $=2$, 
$\dim H^1(\Gamma, \Re^3_{\vect}) = 6 \genus(M) - 6$. 

Let $d^{\nabla}$ be defined on symmetric 2--tensors by $d^{\nabla} h_{abc} =
\nabla_c h_{ab} - \nabla_b h_{ac}$. Thus equation (\ref{eq:codazzi}) can be
written as $d^{\nabla} K = 0$. Tensors satisfying $d^{\nabla} h = 0$ are
called Codazzi tensors. The following proposition shows the relation between
trace--free Codazzi tensors and deformations.
\begin{prop}[Lafontaine \cite{lafontaine:conformal}]
Let $M$ be a compact hyperbolic manifold of dimension $n \geq 2$, with
holonomy representation $\Gamma \subset \SO(n,1)$. Then 
$$
\HH^1 (\Gamma,
\Re^{n+1}_{\vect}) \cong \ker d^{\nabla} \cap \ker \tr .
$$
\end{prop}
In case $n=2$, the trace--free Codazzi tensors are precisely the holomorphic
quadratic differentials. 
Thus, for $n \geq 3$, the deformation space of flat Lorentz structures has
the same dimension as the space of trace--free Codazzi tensors on $M$, 
while for $n
= 2$, the space of flat Lorentz structures has the dimension twice that of
the Teichmuller space of $M$. 

Let $h \in \ker d^{\nabla} \cap \ker \tr$. A cocycle $\alpha \in Z^1(\Gamma,
\Re^{n+1}_{\vect})$  corresponding
to $h$ may be constructed as follows \cite{lafontaine:conformal}. Let $\tilde 
h$ be the lift of $h$ to the universal cover $\HH^n$ of $(M,g_0)$. 
Then $\tilde h = 
\nabla^2 f - f g_0$, for some $f \in C^{\infty}$, cf. 
the note by D. Ferus in \cite{Codazzidisc}. Then $\alpha(\gamma) = f \circ
\gamma - f$ is a cocycle. The function $f$ is unique up to addition by a
function $\beta$ which is the restriction to $\HH^n \subset \Re^{n+1}$ of a
linear function. If $f$ is replaced by $f_1 = f + \beta$, then $\alpha$ is
replaced by $\alpha_1 (\gamma) = \alpha (\gamma) + \beta \circ \gamma -
\beta$, which represents the same class in $H^1(\Gamma, \Re^{n+1}_{\vect})$ as 
$\alpha$, since $\beta \circ \gamma - \beta$ is exact. One shows that
$\alpha$ is nontrivial if $h \ne 0$.

Let $(M,g_0)$ be compact hyperbolic with metric $g_0$ of sectional curvature
$-1$. 
Then the Lorentz cone $(V, \eta)$ over $M$, $V = M \times (0,\infty)$, 
$\eta = - d \rho^2 + \rho^2 g_0$, is a flat, globally
hyperbolic space--time, cf. Proposition \ref{prop:lorentzcone}.  
If $M$ has a two--sided nonsingular totally geodesic hypersurface, then
$H^1(\Gamma, \Re^{n+1}_{\vect}) \ne 0$, cf. \cite{johnson:millson}. Let
$\Sigma \subset M$ be such a totally geodesic hypersurface, and let $h$ be
the induced metric on $M$. Then the Lorentz cone $(V', \eta')$ over $(\Sigma, 
h)$, 
$V' = \Sigma \times (0,\infty)$, $\eta' = -d\rho^2 + \rho^2 h$,  
is totally geodesic in $(V,\eta)$. It follows that replacing $V'$ by 
a Kasner type space--time $W = V' \times (0,R)$, with line element 
$$
- d\rho^2 + \rho^2 h + dr^2,
$$
gives a flat space--time which is a deformation of $(V,\eta)$. 

In dimension $3$, it is still an open question whether or not the space of
flat conformal structures is a manifold. It has been conjectured by Kapovich
\cite{kapovich:deform}, 
that the space of flat conformal structures on a compact hyperbolic
manifold is smooth in dimension 3. 

In the case of flat Lorentz 
structures, on the other hand, the situation is simplified due to the fact
that the Lie algebra structure on $\Re^{n+1}$ is trivial. It follows from
this that the obstructions are all trivial and that the deformation 
space of flat Lorentz
structures is smooth in all dimensions. 

\subsection{Space--times}
We will use both terms ``Lorentz manifold'' and ``space--time''
interchangeably. A space--time $(V,\eta)$ is said to be vacuum if it is Ricci
flat. 

A space--like hypersurface $M$ in a Lorentz manifold $V$ is said to be a
Cauchy surface if $V$ is the domain of dependence of $M$, i.e. if every
inextendible causal curve in $V$ intersects $M$. If $V$ has a Cauchy surface, 
we say that $V$ is globally hyperbolic. Given a vacuum 
space--time $V$, a vacuum 
extension of $V$ is a vacuum 
space--time $V'$ of the same dimension as $V$, 
together with an isometric
imbedding $i: V \to V'$, so that $i(V)$ is a proper subset of $V'$. Inclusion 
induces a partial order on the class of space--times and also on the class of 
globally hyperbolic vacuum space--times, 
so that it is natural to talk about the
maximal globally hyperbolic vacuum extension of a space--time $V$. These
concepts apply directly also to the class of globally hyperbolic flat
space--times. 

If $(M,g,K)$ is a space--like hypersurface in a vacuum space--time $V$, 
with induced Riemann metric $g$ and second fundamental form $K$, then $(g,K)$ 
is a vacuum data set. Similarly,
in case $V$ is flat, we say that the induced data $(g,K)$ is a flat data
set. See section \ref{sec:redhamvol}.
If $V$ is a vacuum (flat) space--time, such that $(M,g,K)$
embeds in $V$ as a Cauchy hypersurface, then $V$ is called a vacuum (flat)
globally hyperbolic extension of $(M,g,K)$. 

In the class of vacuum space--times, the maximal vacuum globally hyperbolic 
extension is unique up to isometry, 
as a consequence of uniqueness for the 
Cauchy problem for the Einstein equations in space--harmonic
coordinates. The following Lemma shows that this holds also in the flat case. 
\begin{lemma}
\begin{enumerate}
\item \label{part:spt} 
Every flat globally hyperbolic
space--time has a unique maximal globally hyperbolic flat extension.
\item \label{part:hypsurf}
Let $(M,g,K)$ be a flat data set. Then there is a maximal globally hyperbolic 
flat extension $V$ which is 
unique up to isometry. In particular, if $\phi: M \to M$ is a diffeomorphism, 
and $g'=\phi^* g, K' = \phi^* K$, then the maximal globally hyperbolic flat
extensions of $(M,g,K)$ and $(M,g',K')$ are isometric. 
\end{enumerate}
\end{lemma}
\begin{proof}
Let $V$ be a flat, globally hyperbolic space--time. Since $V$ is flat it is
also vacuum. 
By a theorem of Choquet--Bruhat and Geroch \cite{ch-b:geroch:glob}, 
every vacuum globally hyperbolic
space--time has a unique maximal globally hyperbolic vacuum extension $V'$. 

Let $W$ denote the Weyl tensor of $V'$. A calculation shows 
$\square W = F(W)$, where $\square$ is the covariant wave operator on
$V'$ and $F$ is a quadratic expression in $W$. 
Uniqueness for the Cauchy problem shows that if $W$ is vanishes on $V$, 
it vanishes on $V'$, and consequently $V'$ is flat. This proves part
\ref{part:spt}. The argument for part \ref{part:hypsurf} is analogous.
\end{proof}
For
brevity we will use the acronym MGHF for maximally globally hyperbolic flat.
In the following
we will let $V$ be a MGHF space--time of
dimension $n+1$, $n \geq 2$, with 
compact Cauchy surface $M$ of hyperbolic type, i.e. $M$ admits a hyperbolic
structure. 

\subsection{Flat space--times}\label{sec:mess}
The presentation in this section follows to a large extent the work
in the unpublished preprint of G. Mess \cite{mess:const:curv},
which dealt with 2+1--dimensional constant curvature space--times. In
particular, most of the proofs in the rest of this section are relatively
straightforward generalizations of the proofs of the analogous results in
\cite{mess:const:curv}, and as the concluding
paragraph of \cite{mess:const:curv} indicates, 
at least some of these generalizations
were known to G. Mess.

A MGHF space--time $V$ with Cauchy surface $M$ has topology $M \times R$, 
and consequently,  $\pi_1 V \cong \pi_1 M$.
A flat Lorentz metric on an $n+1$ dimensional,  
smooth manifold $V$ defines a
$(\ISO(n,1),\Re^{n+1}_1)$--structure on $V$. Letting $\Gamma =
\rho(\pi_1 V)$ be the holonomy representation we can also write $\Gamma =
\rho(\pi_1 M)$. 
The universal cover $\tV$ is
again globally hyperbolic and from this it follows that $\Gamma \subset
\ISO(n,1)_0$, i.e. $V$ is orthochronous. The development $\dev(\tV)$ is a
globally hyperbolic subset of $\Re^{n+1}_1$, $\Gamma$ acts properly
discontinuously on $\dev(\tV)$ and hence on the restriction of the 
$\SO(n,1)_0$ 
bundle over $\Re^{n+1}_1$ to $\dev(\tV)$. It follows that $\Gamma$ is discrete 
in $\ISO(n,1)_0$. 

The following Lemma is a straightforward generalization of \cite[Proposition
4]{mess:const:curv} and \cite[Lemma
2.1]{scannell:rigid}. We include the proof for completeness. 
\begin{lemma}\label{lem:messprop4}
The image of $\dev: \tM \to \Re^{n+1}_1$ embeds $M$ as a 
space--like graph over $\EE^n
= \{(t,x) \in \Re^{n+1}_1 : t = 0 \}$, and $t$ is a proper function on $\tM$. 
\end{lemma}
\begin{proof}
The natural projection $\Re^{n+1}_1 \to \EE^n$ restricts to a
local diffeomorphism $P : \dev(\tM) \to \EE^n$, since
$D(\tM)$ is everywhere space--like. 
As $M$ is compact and strictly 
space--like, the induced Riemann
metric on $\dev(\tM)$ is complete. Further, the pullback to 
$\dev(\tM)$ of the metric
on $\EE^n$ is pointwise larger than the induced metric on $\dev(\tM)$ and is
therefore complete. This shows that 
$P$ is a covering map and hence a global diffeomorphism. Since $\tM$ is
connected, it follows that $\dev(\tM)$ is a graph. 
\end{proof}
An element $\gamma \in \Gamma$ is of the form $x \mapsto Ax + a$. Let $L:
\ISO(n,1)_0 \to \SO(n,1)_0$ be the projection on the linear part, $L(A,a) =
A$, and let $f : \pi_1 M \to \SO(n,1)_0$ be the homomorphism 
$f = L \rho$. Clearly $\ker L = \Re^{n+1}$. 
The short exact sequence 
$$
1 \to \Re^{n+1}_1 \to \ISO(n,1) \overset{L}{\to} \Orth(n,1) \to 1
$$
leads to a short exact sequence 
$$
1 \to T(\Gamma) \to \Gamma \overset{L}{\to} L(\Gamma) \to 1 ,
$$
where $T(\Gamma)$ is the kernel of $L$ restricted to
$\Gamma$. 

For $n=2$, the following Lemma is due to Mess \cite[Proposition
1]{mess:const:curv}. The case $n=3$ follows from \cite{scannell:rigid}. 
\begin{lemma}\label{lem:scannell}
\begin{enumerate}
\item \label{part:isom}
$T(\Gamma) = 0$ and $L : \Gamma \to L(\Gamma)$ is an isomorphism. 
\item \label{part:discret}
$L(\Gamma)$ is a discrete cocompact subgroup of $\SO(n,1)_0$. 
In case $n \geq 3$, $f$ is unique up 
to conjugation.
\end{enumerate}
\end{lemma} 
\begin{proof} 
Since the image under the development map of $\tM$ is
space--like, $T(\Gamma)$ must consist of space--like vectors. 
Since $\Gamma$
is discrete, and $T(\Gamma)$ space--like, 
it follows $T(\Gamma) = Z^k$, some $k \leq n$. By assumption, $\Gamma =
\pi_1 M$ and $M$ is of hyperbolic type. This means $\Gamma$ has no normal
subgroups $Z^k$ and hence $T(\Gamma) = 0$. This proves part \ref{part:isom}.

For part \ref{part:discret}, we use an argument due to Scannell \cite{scannell:rigid},
\cite{scannell:private}. 
Suppose for a contradiction, $L(\Gamma)$ is indiscrete and let $\bLGam$ 
denote the closure of $L(\Gamma)$ in
$\SO(n,1)$. By \cite[Theorem 8.24]{raghunathan:book}, the identity component
$\bLGam_0$ 
is solvable. By \cite[Theorem 5.5.10]{ratcliffe:foundations},
$\bLGam_0$ is elementary, and it follows from the classification of
elementary subgroups in \cite[\S 5.5]{ratcliffe:foundations} that the set of
points $F$ on the sphere at infinity $\partial\HH^n$, left invariant by 
$\bLGam_0$ must consist of one or two points.
The stabilizer in $\SO(n,1)_0$ of a point at infinity is isomorphic to the
group $\Sim^+(\Re^{n-1})$ or orientation preserving similarity
transformations of $\Re^{n-1}$, which is a solvable group. Since $L(\Gamma)$
normalizes $\bLGam_0$, it leaves $F$ invariant and is therefore solvable,
a contradiction to the assumption that $\Gamma = \pi_1 M$ for $M$ of
hyperbolic type. The rigidity statement follows from Mostow rigidity. 
\end{proof}

\begin{prop}\label{prop:lorentzcone}
Given $f: \pi_1 M \to \SO(n,1)_0$, discrete and faithful, there is a
MGHF space--time $V = M \times (0,\infty)$ 
with holonomy $\rho: \pi_1 V \to \ISO(n,1)$ equal to $f$. 
\end{prop}
\begin{proof}
Let $\Gamma = f\pi_1 M$. Then $\Gamma$ acts discontinously on $I^+(\{0\})$,
the interior of future light cone at $0 \in \Re^{n+1}_1$. The quotient $V =
I^+(\{0\})/\Gamma$ is the required space--time. $V$ is globally hyperbolic
since $I^+(\{0\})$ is. 
\end{proof}
The flat space--time $V$ is the Lorentz cone over the compact 
hyperbolic manifold $M$. 
We will refer to such space--times as hyperbolic cone space--times.
The metric on a hyperbolic cone space--time
$V = (0,\infty) \times M$ has the warped product form 
$$
\eta = -d\rho^2 + \rho^2 g_0 ,
$$ 
where $g_0$ is the
hyperbolic metric on $M = \HH^n / \Gamma$. Here $\rho$ corresponds to the
Lorentz distance to $0$ in $\Re^{n+1}_1$. The hypersurface $\rho = s$
has constant mean curvature $\tau = - n/s$, so that $\{ \rho = 1\}$ is 
$\HH^n$. 
\subsection{Existence of strictly convex Cauchy surfaces}
For $\gamma \in \pi_1 M$, $\rho(\gamma) x = f(\gamma) x + t_{\gamma}$,
where $t:\ \pi_1 M \to \Re^{n+1}$, $\gamma \mapsto t_{\gamma}$ is a
1--cocycle, i.e. $t_{\alpha\beta} = t_{\alpha} + f(\alpha) t_{\beta}$. 

The quotient space $H^1(f(\pi_1 M), \Re^{n+1} )$ of the cocycles by the
coboundaries, corresponds to the space of conjugacy classes of
representations $\rho: \pi_1 M \to \ISO(n,1)_0$, such that $L\rho = f$. 
We fix a linear identification $u \mapsto t(u)$, 
of $H^1(f(\pi_1 M), \Re^{n+1} )$ with a subspace of
the space of cocycles $Z^1(f(\pi_1 M), \Re^{n+1})$. 
\begin{prop}[Analogue of \protect{\cite[Proposition 3]{mess:const:curv}}]\label{prop:messprop3}
Given $M$ compact of hyperbolic type and $f: \pi_1 M \to \SO(n,1)_0$, the
holonomy representation for a hyperbolic structure on $M$, let $V$ be the
hyperbolic cone space--time with holonomy $H = f \pi_1 M$. 
Then for any bounded
neighborhood $U$ of $0 \in H^1(H,\Re^{n+1})$, there exists $C> 0$ so that,
letting 
$$
V_C = \{ (t,x) : - t^2 + |x|^2 \leq - C , \quad t > 0 \} /H
$$
there exists a family of flat Lorentz metrics on $V_C$, parametrized by $U$,
so that the space--time $V_C(U)$, corresponding to $u \in U$ has holonomy
$\rho: \pi_1  M \to \ISO(n,1)_0$, sich that $L\rho = f$, and $\rho(\gamma)
x = f(\gamma) x + t_{\gamma}(u)$, where $t_{\gamma}(u)$ is a cocycle
representing $u \in H^1(\pi_1 M, \Re^{n+1})$. We say that $V_C(u)$
represents $u$. Moreover the space--times $V_C(u)$ are future complete and
$M$ is locally strictly convex in each $V_C(u)$. 
\end{prop}
\begin{proof}
Consider the space--time 
$$
V' = \{ (x,t) , x \in \Re^n, t > 0, \quad 
-\half \geq -t^2 + |x|^2 \geq -2 \}/ f \pi_1 M .
$$

Due to the fact that $H^1(f(\pi_1 M), \iso(n,1)) = H^1(f(\pi_1 M),
\Re^{n+1})$, where $\Re^{n+1}$ has the trivial Lie algebra structure (cf. the
discussion above), the obstructions to integrating a ray in
$Z^1(f(\pi_1(M),\Re^{n+1})$ vanish, see \cite{goldman:millson} for
background, and hence by \cite{artin:analytic}, each ray corresponds to a
curve in the deformation space of flat Lorentz structures. It follows from
this that  
there is a neighborhood $U_0$ 
of $0$ in $H^1(\pi_1 M, \Re^{n+1})$, 
and a family of $\ISO(n,1)_0, \Re^{n+1})$ structures on $V'$, parametrized by 
$U_0$, such that the holonomy $\rho: \pi_1 V'(u) \to \ISO(n,1)_0$ represents
$u$. This can also be seen directly by considering the parametrized
development map  $\dev: \tM' \times U_0 \to \Re^{n+1}$, defined by 
$$
\dev(\gamma \cdot x, v ) = f(\gamma) \cdot x + t_{\gamma}(v) .
$$
Choose $\lambda > 0$ so that $U \subset \lambda U_0$ and  
let
$$
V'' = \lambda \cdot V' = \{ (t,x) \in \Re^{n+1}, t > 0, \quad 
-\frac{\lambda^2}{2} \geq - t^2 + |x|^2 \geq - 2\lambda^2 \}/ \pi_1 M .
$$
If we define 
$\dev' : \tV'' \times \lambda \cdot U_0 \to \Re^{n+1}$ by 
$$
\dev'(\lambda\cdot x, \lambda \cdot u) = \lambda \cdot (\dev(x,u)) ,
$$
then 
\begin{align*}
\dev'(\gamma x, u) &= \lambda\dev(\gamma x/\lambda, u/\lambda) \\
&= \lambda(f(\gamma)  x/\lambda + t_{\gamma}(u/\lambda) ) \\
&= f(\gamma) x + t_{\gamma}(u) .
\end{align*}
Thus $\dev'$ is the parametrized developing map of a family of flat
Lorentz structures. We observe that if $U_0$ is chosen sufficiently small,
then for 
$$
N = \{ - t^2 + |x|^2 = - 1 \} / f \pi_1 M , 
$$
$N$ remains locally strictly 
convex in the Lorentz structure determined by any $u \in U_0$. It follows
that $\dev(\tN, u)$ is a strictly convex complete space--like
hypersurface. 
\end{proof}
The following corollary, which is an immediate consequence of the proof of
Proposition \ref{prop:messprop3}, 
provides barriers which will be used to control the
limit of constant mean curvature hypersurfaces in the expanding
direction. 
\begin{cor}\label{cor:barrier}
Let $\eps > 0$ be given. Then, 
in the situation of Proposition \ref{prop:messprop3}, by choosing $U_0$
sufficiently small, the hypersurfaces $N_{s}$, $3/4 \leq s \leq 5/4$,
defined by 
$$
N_s = \{ - t^2 + |x|^2 = - s^2 \} , 
$$
are strictly locally convex in the Lorentz structure corresponding to $u \in
U_0$, and the mean curvature $\tau_{s,u}$ of $N_{s}$ satisfies 
$$
| \tau_{s,u} + n/s | < \eps .
$$
Further $N_{s_1}$ is in the future of $N_{s_2}$ for $s_1 \geq s_2$.
\end{cor}

\begin{definition}[Standard space--time]
A flat globally hyperbolic space--time $V$ with compact hypersurface
$M$ of hyperbolic type is a {\bf standard space--time} if $V$ is the
quotient of the future in $\Re^{n+1}_1$ of a complete space--like strictly
convex hypersurface $\tM$ by a group of Lorentz isometries acting
cocompactly on $\tM$. 
\end{definition}
Given $V$ with compact Cauchy surface $M$, we may assume $t \nearrow \infty$
on $\tM$, after a time orientation. We say $M$ is future oriented. 
\begin{prop}[Analogue of \protect{\cite[Proposition 5]{mess:const:curv}}]
\label{prop:messprop5}
Suppose $V$ is a flat space--time with compact Cauchy surface $M$ of
hyperbolic type. There exists a standard space--time $V'$ containing a
future directed strictly convex surface $M'$ of the same topology and so
that $M$ and $M'$ have the same holonomy representation.
Furthermore, $\tM'$ can be chosen to lie in the future of $\tM$.
\end{prop}
\begin{proof}
By Lemma \ref{lem:scannell} and Proposition \ref{prop:messprop3}, there is a
standard space--time $V'$ containing a locally 
strictly convex hypersurface $M'$ such that, after fixing some identification 
between $\pi_1 M$ and $\pi_1 M'$, the holonomy homomorphisms from $\pi_1 M$
and $\pi_1 M'$ are equal. Replace $M'$ by $M'(K)$, the hypersurface of points 
which lie on the future pointing normals to $M'$ at proper time $K$. $M$ is
compact, so there exists $K$ sufficiently large so that for each $p \in \tM$, 
some point of $\tM'(K)$ lies in the future of $p$. Then $\tM'(K)$ lies
entirely in the future component of $\Re^{n+1}_1 \setminus \tM$, because no
time--like line joins two points in $\tM'$. 
\end{proof}

\begin{prop}[Analogue of \protect{\cite[Proposition 6]{mess:const:curv}}]
\label{prop:messprop6}
Let $V$ be a flat space--time with compact Cauchy surface $M$. Then there is
a flat space--time $V''$, containing a neighborhood of $S$ in $V$, and
containing a standard space--time $V'$. 
\end{prop}

\begin{proof}
Let $\tM'$ be a space--like strictly
convex hypersurface in $V'$ so that for each $m \in \tM$, some point of
$\tM'$ is in the future of $m$ (i.e. $\tM$ is in the past domain of
dependence of $\tM'$) and $\tM'$is also invariant under the holomony
represenation $\rho(\pi_1 M)$ of $V$. It follows from Proposition
\ref{prop:messprop5} 
that there exists such $\tM'$. $\tM'$ together with its
future is the development of our standard space--time $V'$. To construct
$V''$, it is sufficient to show that $\rho(\pi_1 M)$ acts properly
discontinuously on the region $R$ between the disjoint space--like
hypersurfaces $\tM$ and $\tM'$. We then adjoin $R/\rho(\pi_1 M)$ to $M'$ and
thicken the resulting manifold slightly to the past. 

Define a map $h: \tM \times [0,1] \to R$ by $(m,u) \to (1-u) m +
u\varphi(m)$ where $\varphi(m)$ is the intersection of the future pointing
normal to $\tM$ at $m$ with $\tM'$. The map $h$ is proper: a compact subset
$K$ of $R$ lies in a region $t \leq C$ for some $C$, so $h^{-1} K \subset
\tM$ is contained in the set $A = \{ m \in \tM : t \leq C\}$. By Lemma
\ref{lem:messprop4}, $A$ is compact. 

To see that $\rho(\pi_1 M)$ acts properly discontinously on $R$, let $K
\subset R$ be compact. Then for all but finitely many $g \in \rho(\pi_1 M$,
$$
[g \cdot h^{-1}(K)]\cap h^{-1}(K) = \emptyset = h^{-1}(gK) \cap h^{-1}(K) , 
$$
so $gK \cap K = \emptyset$. $\rho(\pi_1 M$ is torsion free so $R$ covers
$R/\rho(\pi_1 S)$. 
\end{proof}
The following Theorem is a direct consequence of Proposition
\ref{prop:messprop6}. 
\begin{thm}\label{thm:strictconv}
Let $V$ be an MGHF space--time with compact Cauchy surface $M$ of hyperbolic
type. Then $V$ contains a locally strictly convex Cauchy surface to the
future of $M$. 
\end{thm} 

\section{Constant mean curvature foliations}

\subsection{Existence of a CMC hypersurface}
A hypersurface $M$ in a space--time $V$ is said to have constant mean
curvature (CMC) is $\tr K$ is constant on $M$. CMC hypersurfaces 
satisfy a space--time maximum principle, and existence of CMC hypersurfaces
can be proved using barriers. We use the existence of a strictly convex
hypersurface to construct barriers.  
\begin{thm}\label{thm:cmcexist}
Let $V$ be a MGHF space--time with compact Cauchy surface $M$ of
hyperbolic type. Then $V$ contains a CMC Cauchy surface to the future of $M$. 
\end{thm}
\begin{proof}
By Theorem \ref{thm:strictconv}, $V$ contains a locally strictly convex
Cauchy surface $M'$ to the future of $M$. 
Let $K^a_{\ b}(t)$ be the mixed for of the second fundamental form of the Gauss
foliation $\{S_t\}_{t=0}^{\infty}$, to the future 
(the expanding direction) of $M'$,
with the sign convention 
$$
\partial_t g_{ab} = - 2 K_{ab} , 
$$
where $g_{ab}$ is the induced Riemann metric on $S_t$. $K^a_{\ b}$ is the
shape operator of the Gauss foliation and satisfies the Riccati equation 
$$
\partial_t K^a_{\ b} = K^a_{\ c} K^c_{\ b}  .
$$
By diagonalizing this is seen to have the solution 
\begin{equation}\label{eq:riccati-sol}
K(t) = (K^{-1}(0) - t)^{-1} , 
\end{equation}
and it follows that $\tr K(t)$ is strictly increasing with $t$, and $\tr K(t)
\nearrow 0$ as $t \to \infty$. Since $\tr K(0) < C < 0$, it follows that
there is a $t_0$ so that $M', S_{t_0}$ form barriers for the CMC problem and
by \cite{gerhardt:CMC}, there is a CMC Cauchy surface $M''$ in $V$ to the future
of $M$. The mean curvature $\tau$ of $M''$ can be choosen to have any value
in the open interval $(\max \tr K(0), 0)$. 
\end{proof}

\subsection{CMC Cauchy surfaces are strictly convex}
By results of Treibergs and Choi, entire CMC hypersurfaces in $\Re^{n+1}_1$
have the property that the second fundamental form $K$ is semidefinite, and
if $K$ has a zero eigenvalue at some point, then the CMC hypersurface
splits. We use this to prove the following Proposition.
\begin{prop}\label{prop:CMCconv}
Let $V$ be a MGHF space--time with compact CMC Cauchy surface $M$ of
hyperbolic type. Then $M$ is strictly locally convex.  
\end{prop}
\begin{proof}
We may assume the mean curvature $\tau$ of $M$ satisfies $\tau <
0$. Consider $\tM$, the development of $M$ in $\Re^{n+1}$. By
\cite{treibergs:cmc}, the second fundamental form $K$ of $\tM$ is negative
semi--definite. We must prove that $K$ is negative definite. 
Suppose for a contradiction, there exists $m \in \tM$ such that $K$ has a
zero eigenvalue. Then by \cite{treibergs:choi}, 
$\tM$ splits metrically, 
$\tM = \tM^{n-k} \times \Re^k$ for some $k$, $1 \leq k \le n$, where
$\tM^{n-k}$ has negative definite second fundamental form and the same mean
curvature as $\tM$. In particular, the rank of $K$ is constant. 

The tangent spaces split $T_m \tM = E_{1,m} \oplus E_{2,m}$,
where $E_{1,m}, E_{2,m}$ are of dimension $n-k$ and $k$, respectively. 
If we let $P$ denote the orthogonal projection on $\Re^k$, $E_{2,m} = P T_m
\tM$. According to
this splitting, $K = K_1 \oplus 0$, where $K_1 < 0$. 

Let 
$\rho : \pi_1 M \to \ISO(n,1)_0$ be the holonomy
representation of $V$. For $\gamma \in \pi_1 M$, $\rho(\gamma)$ is an
affine isometry of $\Re^{n+1}_1$, $\rho(\gamma) m = f(\gamma) m +
t_{\gamma}$, 
where $f = L \rho$. 

Now we show that $f(\gamma)$ commutes with $P$, i.e. $f(\gamma)
E_{1,m} = f(\gamma) E_{1, \rho(\gamma) m}$, $f(\gamma) E_{2,m} =
E_{2,\rho(\gamma) m}$. 
To see this, use the fact that $\rho(\gamma)^* K  = K$ by construction. This
implies that if we choose $v \in E_{2,m}$, $w \in T_m\tM$, 
$$
0 = \rho(\gamma)^* K_m(v,w) = K_{\rho(\gamma) m}(\rho(\gamma)_* v,
\rho(\gamma)_* w ) =  K_{\rho(\gamma) m}(f(\gamma) v, f(\gamma) w ) .
$$
But $f(\gamma): T_m \tM \to T_{\rho(\gamma) m} \tM$ is an isomorphism, so
$f(\gamma) E_{1,m} = E_{1,\rho(\gamma) m}$. Similarly, one sees that 
$f(\gamma) E_{2,m} = E_{2,\rho(\gamma) m}$.

The projection of $f$, $\bar f = Pf$ is
a homomorphism $\pi_1 M \to \SO(k)$. But $f(\pi_1 M)$ is the holonomy of
a hyperbolic structure in dimension $n$. 
This gives a contradiction, and we
can conclude that $M$ is locally strictly convex. 
\end{proof}
A CMC foliation $\FF$ of a space--time $V$ spatially compact Cauchy surface,
is a smooth foliation 
$\FF = \{M_{\tau} : \tau \in I\}$ in $V$, whose leaves are CMC Cauchy
surfaces. Let $V_{\FF}$ be the maximal domain
of $V$ foliated by $\FF$. We say that $V$ has a global CMC foliation if
$V_{\FF} = V$. 
As a consequence of the strict local convexity of the CMC Cauchy surface, we
get a local foliation by CMC Cauchy surfaces.
\begin{cor}\label{cor:CMClocal}
Let $V$ be a MGHF space--time, with a compact CMC Cauchy surface $M$ of
hyperbolic type, with mean curvature $\tau < 0$. Then there is an $\eps >
0$, and a local CMC
foliation $\FF = \{M_t : t \in (\tau-\eps, \tau+\eps) \}$. 
\end{cor}
\begin{proof} In view of the strict local convexity of $M = M_{\tau}$, the
Gauss foliation w.r.t. $M_{\tau}$ can be extended both to the past and the
future. Using the solution (\ref{eq:riccati-sol}) to the Riccati equation,
we see that there is an $\eps > 0$ so that the leaves in the Gauss foliation
form barriers with $\tr K < \tau - \eps$ to the past, and $\tr K > \tau +
\eps$ to the future. Hence by \cite{gerhardt:CMC}, we get a local foliation
by CMC hypersurfaces with mean curvature taking all values in the interval
$(\tau-\eps, \tau+\eps)$. 
\end{proof}

\subsection{Global CMC foliation}\label{sec:cmcglobal}
Next we show that the local CMC foliation can be extended to a global
one. This follows from \cite{andersson:etal:2+1grav} in case $n=2$, and 
\cite[Theorem 0.1]{anderson:foli} in case $n=3$. For the general case, we
will use the results of \cite{andersson:moncrief:ellhyp}. The method of
\cite{anderson:foli} may also be adapted to this situation, with some work. 
\begin{thm}\label{thm:cmcglobal}
Let $V$ be a MGHF space--time with compact Cauchy surface $M$ of hyperbolic
type. Then $V$ has a global CMC foliation. 
\end{thm}
\begin{proof}
Since $V$ is flat, the Einstein CMC evolution equations with zero shift, in
CMC time $t = \tr K$, 
take the form 
\begin{subequations}\label{eq:evol}
\begin{align}
\partial_t g_{ab} &= - 2 N K_{ab} , \label{eq:g-evol} \\
\partial_t K_{ab} &= - \nabla_a \nabla_b N - N  K_{ac} K^c_{\ b} ,  \\
- \Delta N + |K|^2 N &= 1  .
\end{align}
\end{subequations}
Let $\UU = (u,v)$, $u_{ab} = g_{ab}$, $v_{ab} = - 2 K_{ab}$, 
and let 
\begin{align}
L \UU &= \begin{pmatrix} \partial_t u - N v \\
                        \partial_t v \end{pmatrix} ,  \\
F_1[u,v] &= 0 , \\
F_2[u,v] &= 2 \nabla_a \nabla_b N + \half N v_{ac} v^c_{\ b} ,
\end{align}
where $\nabla$ is defined w.r.t. $u$, and indices are raised and lowered
w.r.t. $u$. Let $F = (F_1, F_2)$. 
Then (\ref{eq:evol}) is an elliptic--hyperbolic system of the form 
$L\UU = F$, considered in \cite{andersson:moncrief:ellhyp}. Note that
$L$ is not hyperbolic in the usual sense, since the principal part is 
$\partial_t$. However, all the basic energy estimates of
\cite{andersson:moncrief:ellhyp} are valid for this system, and consequently,
the local existence and continuation properties hold for this system. The
relevant continuation principle can be stated as follows. Let 
$(T_-, T_+)$
be the maximal interval such that there is a CMC foliation $\FF$ with mean
curvatures taking all values in $(T_-, T_+)$. Then either $(T_-,T_+) =
(-\infty, 0)$ or the quantity 
$$
\max(\Lambda[g], |g|_{C^1}, |N|_{C^1}, |N^{-1}|_{L^{\infty}}, |N|_{C^1})
$$
diverges, as $t \to T_-$ or $t \to T_+$. 
Here $\Lambda[g]$ denotes the ellipticity constant of $g$ and
$|\cdot|_{C^1}$ denotes the $C^1$ norm.   

Let $\tau = \tr K$. 
By the work of Treibergs \cite{treibergs:cmc}, the estimate
\begin{equation}\label{eq:Kbound}
|K|^2 \leq \tau^2 
\end{equation}
holds. It follows that $|\Ric| \leq C_n \tau^4$. 
From the bound on $|K|$ we get using the maximum principle and the Lapse
equation 
$1/\tau^2 \leq N \leq n/\tau^2$. It follows from this and the evolution
equation for $g$ that the ellipticity constant and the quantities
$|N|_{L^{\infty}}, |N^{-1}|_{L^{\infty}}$ can't diverge as long as $\tau \in
(-\infty, 0)$. It remains to consider $|g|_{C^1}$ and $|N|_{C^1}$. We will
use the results of Anderson and Cheeger \cite{anderson:cheeger}. In view
of the $|K|$ bound and the evolution equation for $g$, collapse is ruled
out for $\tau \in (0,\infty)$. 
By the Ricci bound $|\Ric|\leq C_n \tau^4$, we have $g$ under control in
$C^{1,\alpha}$ in harmonic coordinates and the harmonic coordinate radius
$r_h$ is bounded away from zero. This takes care of $|g|_{C^1}$. In view of
this regularity of $g$, and 
standard elliptic estimates we get from the Lapse equation that 
$|N|_{W^{2,p}}$  is bounded for any $p > 1$, and hence $|N|_{C^1}$ is
bounded. This proves in view of the continuation principle that the CMC
foliation extends to all values of $\tau$ in $(-\infty,0)$. In view of the
Lapse bounds, it follows that proper time runs to infinity along the
foliation in the expanding direction $\tau \nearrow 0$, and in view of 
the fact that focal distance tends to zero as $\tau \searrow -\infty$, the
foliation covers all of $V$. 
\end{proof}

\section{The limit in the expanding direction}
In this section, we consider the limit of $M_{\tau}$ as $\tau \to 0$. 
We will consider the effect of scaling
$V$ and $M$. 

\begin{thm}\label{thm:limit}
Let $(V,\eta)$ be a MGHF space--time with Cauchy surface $M$ 
of hyperbolic type.
Let $M_{\tau}$ be the Cauchy surface with mean curvature $\tau < 0$, and let
$g(\tau)$ be the induced metric on $M_{\tau}$. Then as
$\tau \nearrow 0$, $(M,\frac{\tau^2}{n^2} g(\tau))$ converges in $C^0$ to 
$(\HH^n , g)$, where $g$ is the hyperbolic metric with sectional curvature
$-1$. Further the rescaled volume 
$\frac{|\tau|^n}{n^n} \Vol(M,g_{\tau})$ tends to the hyperbolic
volume, i.e. 
$$
\lim_{\tau \nearrow 0}\frac{|\tau|^n}{n^n} \Vol(M,g_{\tau}) \to \Vol(M,g_0) .
$$
\end{thm}
\begin{proof}
Let $U_0 \subset H^1(\pi_1 M, \Re^{n+1})$ be a small neighborhood of zero as
in the proof of Proposition \ref{prop:messprop3}, and let $u \mapsto t(u)$ be 
the map from $U_0$ to the space of cocycles, with corresponding parametrized
development map $(\gamma , u ) \mapsto \dev(\gamma, u)$, with corresponding
parametrized holonomy representation $(\gamma, u) \mapsto \rho(\gamma ,
u)$. By the proof of Proposition \ref{prop:messprop3}, there is a $v$ such
that the development map of $V$ is given by 
$$
\dev(\gamma, v) x = f(\gamma) x + t_{\gamma}(v) , 
$$
with cocycle
$t(v)$. Then there exists  $\lambda >
0$, $u \in U_0$, such that $\dev( \cdot , \lambda \cdot u) = \dev( \cdot ,
v )$ is the development map of $V$ with cocycle $t(v) = \lambda t(u)$. 

Instead of scaling up by a large $\lambda$ as in the proof of Proposition
\ref{prop:messprop3}, we will here scale down by $\lambda^{-1}$, for large
$\lambda$. 

We denote by $V_{\lambda}$, the space--time $V$ scaled by a factor
$\lambda^{-2}$. 
The effect of the scaling of $V$ by $\lambda^{-2}$ 
can be described in two ways. First, the
scaling of the cocycle is given by $t(v) \to t(\lambda^{-2} v) =
\lambda^{-2} t(v)$, so as $\lambda \to \infty$, the geometric structure 
of $V_{\lambda}$ converges to the space--time $V_{\infty}$ with holonomy
representation $f$. Secondly, the scaling can be described by $\eta \to
\eta_{\lambda} = 
\lambda^{-2} \eta$, a scaling of the space--time metric. 

Let $M_{\tau}$ be a CMC Cauchy surface of $V$ with mean curvature $\tau$, 
with induced metric and second
fundamental form $(g,K)$. The  scaling of $V$ by $\lambda^{-2}$ scales the
data on $M$ by 
$g \to \lambda^{-2} g$, $K \to \lambda^{-1} K$, $\tau \to \lambda
\tau$. Therefore, if we consider a sequence of CMC Cauchy surfaces 
with mean curvatures $\tau_i$, such that
$\tau_i \nearrow 0$, scaling $V$ by 
$\tau_i^2/n^2$ has the effect of sending $M_{\tau_i}$ to a 
CMC Cauchy surface $M_i$, with mean curvature $-n$,  
in the space--time $V_i = V_{\frac{n^2}{\tau_i^2}}$. 

The sequence of space--times $V_i$
converges to the space--time $V_{\infty}$, as $i \to \infty$. For $i$
sufficiently large, we may apply Corollary \ref{cor:barrier} to find CMC
hypersurfaces $N_{1,i}, N_{2,i}$ with mean curvatures $\tau_{1,i},
\tau_{2,i}$, satisfying
$$
\tau_{1,i} < -\frac{4n}{3} + \eps , \qquad \tau_{2,i} > - \frac{4n}{5} - \eps ,
$$ 
for some small $\eps > 0$. 
The hypersurfaces 
$N_{1,i}, N_{2,i}$ correspond to $s = 3/4, 5/4$, respectively,
in Corollary \ref{cor:barrier}. Further, $N_{1,i}$ is in the past of
$N_{2,i}$. 
These are therefore barriers for the CMC problem for the mean curvature
$-n$. It follows by \cite{gerhardt:CMC} that there is a CMC Cauchy
surface $\bar M'_i$ with mean curvature $-n$, such that $M'_i$ is in the future
of $N_{1,i}$ and in the past of $N_{2,i}$. By uniqueness of CMC
hypersurfaces 
in space--times with compact
Cauchy surface, cf. \cite{brill:flaherty},  \cite{marsden:tipler}, we see
that $M'_i = M_i$. 
Thus, $M_i$ is in the future of $N_{1,i}$ and in the past of
$N_{2,i}$. 

Now we pass to the universal cover. Let $\tN_{1,i}, \tN_{2,i},
\tM_i$ be the developements under the parametrized development map of
$N_{1,i}, N_{2,i}, M_i$, respectively. By construction, the same causal
relations hold as before.  Since $V_i$ tends to the space--time $V_{\infty}$,
as $i \to \infty$, we find that 
$\tN_{1,i}, \tN_{2,i}$ tend to the hyperboloids with mean curvature $-4n/3,
-5n/4$, respectively, as $i \to \infty$. 
The convergence is uniform on compacts.

It remains to check that the volume of the induced metric tends to the volume
of the hyperbolic metric on the hyperboloid $\HH^n /f \pi_1 M$. 
Recall that by Lemma \ref{lem:messprop4}, $\tM_i$ 
is the graph of a function $\phi_i$ over $\EE^n$.
By the work of Treibergs \cite[p. 54]{treibergs:cmc}, the height function
$\phi_i$ of
$\tM_i$ has uniform $C^3$ bounds on compact subsets, and further 
$|D\phi| < 1$ uniformly on compacts. 
These bounds, applied to $\tM_i$, together with
the above barrier construction, we find that as $i
\to \infty$, $M_i$
converges uniformly on compacts to the hyperboloid with mean curvature $-n$.

A calculation shows that the induced
volume element $\nu_i$ on $\tM_i$ is given in terms of the gradient 
$\nabla\phi_i$ by 
$\nu_i = \sqrt{1 - |\nabla\phi_i|^2}dx^1\cdots dx^n$, 
where $|\nabla \phi_i|$ is
the Euclidean norm of the gradient of $\phi_i$. The $C^3$ bound shows in
view of Arzela-Ascoli that 
$\nabla\phi_i$ converges uniformly on compacts, and hence the volumes 
$\Vol(M_i)$ satisfy 
$$
\lim_{i \to \infty} \Vol(M_i) = \Vol(M, g_0) .
$$
Since the induced metric on $M_i$ is the rescaled metric
$\frac{\tau_i^2}{n^2} g(\tau_i)$, this means precisely that the rescaled
volume converges. 
\end{proof}
\section{The Gauss map}
Let $V$ be a MGHF space--time with compact Cauchy surface $M$ of hyperbolic
type. Let $\norm$ be the future directed
time--like unit normal, $\la \norm , \norm \ra = -1$. Let as usual $\rho:
\pi_1 M \to \ISO(n,1)_0$ be the holonomy representation of $V$ and let $f = L 
\rho$ be the holonomy representation of $M$. 
\begin{prop}
The map $x \mapsto \norm(x)$ defines a well--defined map 
$\varphi: M \to \HH^n/f\pi_1(M)$, the Gauss map. 
\end{prop}
\begin{proof}
Consider $\dev(\tM) \subset \Re^{n+1}_1$. We denote the unit future directed 
time--like normal of
$\dev(\tM)$ by $\tnorm$. Since $\Re^{n+1}_1$ has trivial holonomy, $\tnorm(x)$
may be identified with a point $\tvarphi(x) \in \HH^n$. This 
defines the Gauss map $\tvarphi: \dev(\tM) \to \HH^n$. In view of the fact that 
$\tnorm (\rho(\gamma) x) = \rho(\gamma)_* \tnorm(x) = f(\gamma) \tnorm(x)$, 
the Gauss map has the following
equivariance property, 
$$
\tvarphi(\rho(\gamma) x) = f(\gamma) \tvarphi(x) ,
$$
where in the RHS, we use the action of $f\pi_1 M$ on $\HH^n$. It follows from 
this equivariance property, that the map $\tvarphi$ induces a map 
$\varphi: M \to \HH^n / f \pi_1 M$. This is the required Gauss map. 
\end{proof}

\begin{lemma} Let $M$ be a CMC Cauchy surface in $V$, with Gauss map
$\varphi : M \to \HH^n / f \pi_1 M$. Then $\varphi$ is harmonic. 
\end{lemma}
\begin{proof}
Consider $\tvarphi: \dev(\tM) \to \HH^n$. $\dev(\tM)$ is a CMC hypersurface in 
$\Re^{n+1}_1$ so by \cite{treibergs:choi} $\tvarphi$ is harmonic. The
corresponding property for $\varphi$ follows immediately. 
\end{proof}

\begin{lemma} $\varphi$ is a diffeomorphism, which is isotopic to the
identity. 
\end{lemma}
\begin{proof}
We may without loss of generality assume the mean curvature of $M$ is $-n$.
Clearly $\varphi$ is smooth. Let $t$ be the cocycle so that $\dev(\gamma) x = 
f(\gamma) x + t_{\gamma}$. By scaling the cocycle as in the proof of
Proposition \ref{prop:messprop3}, we get a one parameter family of
space--times and CMC hypersurfaces and Gauss maps, starting at 
$\HH^n/f\pi_1 M, \id$ and ending
at $M,\varphi$. This shows that $\varphi$ is homotopic to the identity.

Let $N$ denote $M$ with its hyperbolic metric $g_0$ and fix compatible
orientations on $M, N$. 
Let $\deg\varphi$ denote the degree of the map $\varphi : M \to N$. Since
$\deg \varphi$
is 
a homotopy invariant, we may in view of the fact that $\varphi$ is homotopic to
the identity, conclude that $\deg\varphi = 1$. 

Recall that $d\varphi: TM \to TN$ is given by the second fundamental form,
$d\varphi = -K$. By Proposition \ref{prop:CMCconv}, $-K$ is positive
definite. It follows from the inverse function theorem that $\varphi$ is
locally smoothly invertible. 
Since $\deg\varphi = 1$, and $d\varphi$ is positive definite, it follows that
$\varphi$ is a bijection, and hence it is a diffeomorphism.  
\end{proof}

The harmonic map energy of $\varphi$, defined by 
$E(\varphi) = \int_M |d\varphi|^2\mu_g$,
can be written as 
$$
E(\varphi) 
= \int_M |K|^2 d\mu_g = \int_M (R + (\tr K)^2 )\mu_g = \int_M R \mu_g 
+ \tau^2 \Vol(M,g) , 
$$
where we used the constraint equation $R + (\tr K)^2 = |K|^2$ with $R$ the
scalar curvature of $(M,g)$, in the last
step. In case $n=2$, $\int_M R \mu_g = 4\pi \chi(M)$ by Gauss-Bonnet, which
gives the interesting formula 
$$
E(\varphi) = 4 \pi \chi(M) + \tau^2 \Vol(M,g) ,
$$
found by Puzio \cite{puzio:gauss}.
In case 
$n \geq 3$, we note $\int_M R \mu_g < 0$, in view of the fact that $M$ is of
Yamabe type $-1$. 

\section{The rescaled volume and the $\sigma$--constant}
Let $V$ be a MGHF space--time with Cauchy surface $M$ of hyperbolic type, and 
let $\FF = \{M_{\tau} : \tau \in (-\infty, 0)\}$ be the global CMC foliation of $V$
as in Theorem \ref{thm:cmcglobal}. In the rest of this section, $M$ will
denote a CMC Cauchy surface, with induced metric $g$ and second fundamental
form $K$, with $\tr K = \tau$, constant on $M$. 
Let 
$$
\Ham(\tau) = \Ham(g,K) = |\tau|^n \Vol (M,g) = |\tau|^n \int_M \mu_g
$$
Then $\Ham(g,K)$ can be thought of as a function on the space of CMC flat
data sets on $M$. 

The infimum $\inf_{g,K} \Ham(g,K)$ of $\Ham$ over all CMC flat data on $M$,
is a quantity that is analogous to the $\sigma$--constant of $M$, but in the
restricted setting of flat CMC data.

Via the Yamabe equation, we may write the metric $g$ as 
$g = u^{\frac{4}{n-2}} h$,
where $h$ has constant negative scalar curvature $-n(n-1)$. Let 
$$
K = \frac{\tau}{n} g + u^{-2}\sigma
$$
where $\sigma$ is satisfies $\tr_h \sigma = 0$, $\nabla^a[h] \sigma_{ab} =
0$. Then $K$ is divergence free w.r.t. $g$, and using the Yamabe equation,
the vacuum 
constraint equation $R_g + (\tr_g K)^2 - |K|_g^2 = 0$, takes the form 
$$
- \frac{4(n-1)}{n-2} \Delta u + R_h u + \frac{n-1}{n} \tau^2
u^{\frac{n+2}{n-2}} - u^{\frac{2-3n}{n-2}} |\sigma|_h^2
$$
Thus we may write 
$$
\Ham(g,K) = |\tau|^n \int_M (u^{\frac{4}{n-2}})^{n/2} \mu_h
$$
where via the constraint equation, we have $u = u(h,\sigma,\tau)$. A maximum
principle argument shows that $u(h,\sigma,\tau) \geq u(h,0,\tau)$, and hence 
we get that $u(h,0,\tau)$ satisfies 
$$
u^{\frac{4}{n-2}} = \frac{n^2}{\tau^2}
$$ 
This gives 
$$
\Ham \geq n^n \Vol(M,h) 
$$
where $h$ is the unique metric in the conformal class of $g$ with scalar
curvature $R_h = - n(n-1)$ and shows that 
$$
\inf \Ham(g,K) = n^n \inf_{h \in \MM_{-n(n-1)}} \Vol(M,h) 
$$ 
where the infimum is over the space of CMC vacuum data on $M$. It follows
that 
$$
\sigma(M) = -\frac{n-1}{n} \inf \Ham^{2/n}
$$
These relations have been studied in case $n=3$ by Fischer and Moncrief
\cite{fischer:moncrief:sigma}.

\begin{thm} Let $M$ be a compact 
manifold of hyperbolic type, of dimension $n \geq 3$, 
and let $g_0$ be the hyperbolic metric on $M$ with sectional curvature $-1$. 

Let $(g,K)$ be flat data set on $M$ with mean curvature $\tau = \tr K$. 
Then the rescaled volume $\Ham(g,K) = |\tau|^n \Vol(M,g)$ satisfies
\begin{equation}\label{eq:Hamineq:thm}
\Ham(g,K) \geq n^n \Vol(M,g_0) .
\end{equation}
Equality in (\ref{eq:Hamineq:thm}) holds if and only if $g$ is isometric to
$\frac{n^2}{\tau^2} g_0$ 
and $K = \frac{\tr K}{n} g$. 
\end{thm}
\begin{proof} We will use a monitonicity property for $\Ham$, which was first 
noticed by A. Rendall (unpublished) and studied in detail by Fisher and
Moncrief \cite{fischer:moncrief:sigma} for dimension $n=3$. 

Let $\hK$ denote the tracefree part of $K$, $K = \frac{\tau}{n} g +
\hK$. By integrating the Lapse equation we see 
$$
\int_M (\frac{d\tau}{dt} - N \frac{\tau^2}{n})\mu_g  = 
\int_M N |\hK|^2 \mu_g
$$
Time differentiation $\Ham(\tau)$ using the evolution equation
(\ref{eq:g-evol}) and shows
\begin{equation}\label{eq:Hamineq}
\partial_t \Ham(\tau) = - n |\tau|^{n-1} \int_M N |\hK|^2 \mu_g \leq 0
\end{equation}
and $\partial_t H(\tau) = 0$ if and only if $\hK = 0$. 

In case $\hK = 0$ at $\tau$, a calculation shows 
$\partial_t^2 \Ham(\tau) = 0$. Further, in this case, 
$$
\partial_t^3 H (\tau) = - |\tau|^{n-1} n \int N |\partial_t \hK|^2 \mu_g
$$
In case $\hK = 0$ at $\tau$, $N = n/\tau^2$ 
by the Lapse equation with $d\tau/dt =1$ and the evolution equations 
(\ref{eq:evol}), with we get 
$$
\partial_t \hK_{ab} = \frac{2}{n} g 
$$
This shows that in case $\partial_t \Ham(\tau) = 0$, then $\partial_t^2
\Ham(\tau) =0$, and 
$$
\partial_3 \Ham(\tau) = - 4n |\tau|^{n-3} \Vol(M,g) < 0
$$
This shows that $\Ham$ is monotonically decreasing. As $\Ham$ is nonnegative, 
it follows that the limit $\lim_{\tau \to 0} \Ham(\tau)$ exists, and 
$$
\lim_{\tau \to 0 } \Ham(\tau) = \inf_{\tau \in (-\infty, 0)} \Ham(\tau)
$$

On the other hand, we have seen that $\tau \to \Ham(\tau)$ is monotone
decreasing, and by Theorem \ref{thm:limit}, 
$\lim_{\tau \nearrow 0} \Ham(\tau) = n^n\Vol(M,g_0)$, where $g_0$ is the
hyperbolic metric on $M$ with sectional curvature $-1$. This proves
(\ref{eq:Hamineq:thm}). 

It remains to prove the rigidity statement. Suppose that $(g,K)$ are flat CMC
data such that $\Ham(\tau) = \Ham(g,K) = n^n\Vol(M,g_0)$. Then $\partial_t
\Ham(\tau) = 0$, which implies that $K = \frac{\tr K}{n} g$. 
After scaling, we may assume
$\tr K = - n$. It follows that 
$$
R_{ab} = -(n-1) g_{ab}
$$
In case $n=3$ this is enough to complete the proof. 
For the higher dimensional case, However, it follows from $\Ham(g,K) = n^n\Vol(M,g_0)$ that $\Vol(M,g) =
\Vol(M,g_0)$, and therefore, by \cite{BCG2}, $(M,g)$ is isometric to $(M,g_0)$.
This completes the proof. 
\end{proof}


\providecommand{\bysame}{\leavevmode\hbox to3em{\hrulefill}\thinspace}

\end{document}